\documentclass[pdflatex,sn-mathphys-num]{sn-jnl}

\usepackage{amsmath,amssymb,amsfonts,amsthm}%
\usepackage{color,dsfont}
\usepackage{comment}
\allowdisplaybreaks

\theoremstyle{thmstyleone}%
\newtheorem{theorem}{Theorem}[section]
\newtheorem{lemma}[theorem]{Lemma}%
 
\theoremstyle{thmstyletwo}%
\theoremstyle{thmstylethree}%
\numberwithin{equation}{section}

\newcommand{\B}{\mathcal{B}}
\newcommand{\K}{\mathcal{K}}
\newcommand{\M}{\mathcal{M}}
\newcommand{\X}{\mathcal{X}}

\newcommand{\N}{\mathbb{N}}
\newcommand{\Z}{\mathbb{Z}}
\newcommand{\R}{\mathbb{R}}
\newcommand{\C}{\mathbb{C}}
\newcommand{\T}{\mathbb{T}}

\newcommand{\Ind}{\operatorname{Ind}}
\newcommand{\wind}{\operatorname{wind}}
\newcommand{\Ker}{\operatorname{Ker}}

\newcommand{\Image}{\operatorname{Im}}
\newcommand{\clos}{\operatorname{clos}}
\newcommand{\Lspan}{\operatorname{span}}
\newcommand{\Log}{\operatorname{Log}}
\newcommand{\Arg}{\operatorname{Arg}}

\newcommand{\essinf}{\operatornamewithlimits{ess\,inf}}
\newcommand{\esssup}{\operatornamewithlimits{ess\,sup}}

\newcommand\restr[2]{{% we make the whole thing an ordinary symbol
  \left.\kern-\nulldelimiterspace % automatically resize the bar with \right
  #1 % the function
  \vphantom{\big|} % pretend it's a little taller at normal size
  \right|_{#2} % this is the delimiter
  }}

\begin{document}

%% Article title
\title{Fredholm criteria for Wiener-Hopf operators with continuous symbols acting on some Banach function spaces}

%% Author name
\author*[1]{\fnm{M\'arcio} \sur{Valente}}
\email{mac.valente@campus.fct.unl.pt}

%% Author affiliation
\affil*[1]{
\orgdiv{Centro de Matem\'atica e Aplica\c{c}\~oes, Departamento de Matem\'atica}, 
\orgname{Faculdade de Ci\^encias e Tecnologia, Universidade Nova de Lisboa}, 
\orgaddress{
\street{Quinta da Torre}, 
\city{Caparica}, 
\postcode{2829--516}, 
\country{Portugal}}
\\ \\
To Professor Sergei Grudsky on the occasion of his 70th birthday
}

%% Abstract
\abstract{\unboldmath Let $X(\R_+)$ be one of the following three Banach function spaces: a Lorentz space $L^{p,q}(\R_+)$ with $1 < p,q < \infty$; a reflexive Orlicz space $L^{\Phi}(\R_+)$; or a variable Lebesgue space $L^{p(\cdot)}(\R_{+})$ with variable exponent $p(\cdot)\in \B_{M}(\R)$. We extend the Fredholm criteria for Wiener-Hopf operators with continuous symbols on the Lebesgue space $L^{p}(\R_{+})$, $1 < p < \infty$, obtained by Roland Duduchava in the late 1970s, to the space $X(\R_{+})$.}

%% Keywords


\maketitle

%%%%%%%%%%%%%%%%%%%%%%%%%%%%%%%%%%%%%%%%%%%%%%%%%%%%%%%%%%%%%%%%%%%%%%%%%%%%
%%%%%%%%%%%%%%%%%%%%%%%%%%%%%%%%%%%%%%%%%%%%%%%%%%%%%%%%%%%%%%%%%%%%%%%%%%%%
%%%%%%%%%%%%%%%%%%%%%%%%%%%%%%%%%%%%%%%%%%%%%%%%%%%%%%%%%%%%%%%%%%%%%%%%%%%%

%% MAIN TEXT

\section{Introduction and main result}\label{sec: introduction}

Let $\X$ be a Banach space, $\B(\X)$ be the Banach algebra of all bounded linear operators on $\X$ and $\K(\X)$ be the closed two-sided ideal of all compact operators in $\B(\X)$. Formally, we say that an operator $T\in \B(\X)$ is Fredholm on $\X$ if its image $\Image T$ is closed in $\X$ and the numbers
\[
n(T) := \dim \Ker T, 
\qquad 
d(T) := \dim \X / \Image T, 
\]
are finite. It is well-known (see, e.g., \cite[Chap.~XI, Corollary~2.3]{GGK90}) that the condition on $\Image T$ being closed is redundant. To each Fredholm operator $T$, we can also assign a constant, called the Fredholm index of $T$, which is the number $\Ind T := n(T) - d(T)$. In practice, this definition is of rather uncomfortable use. Fortunately, an alternative characterization is available (see, e.g., \cite[Chap.~XI, Theorem~5.1]{GGK90}) and affirms that an operator $T$ is Fredholm on $\X$ if and only if there is some $S\in \B(\X)$ such that the operators $I - TS$ and $I - ST$ are compact on $\X$. This classical result ends up providing us with a wonderful revelation: an operator $T$ is Fredholm on $\X$ if and only if it is invertible modulo compact operators.

%%%%%%%%%%%%%%%%%%%%%%%%%%%%%%%%%%%%%%%%%%%%%%%%%%%%%%%%%%%%%%%%%%%%%%%%%%%%

Having introduced the concept that we are going to study, we will now specify the context in which we are going to work, namely, the function spaces and later the operators. First and foremost, given a domain $\Omega\in \{\R_{+}, \R\}$, where $\R_{+} := (0, \infty)$ and $\R$ is the real line, the function spaces that we will consider will be one of the following three Banach function spaces $X(\Omega)$: 
\begin{itemize}
\item[(a)] Lorentz spaces $L^{p, q}(\Omega)$ with $1 < p, q < \infty$;
\item[(b)] reflexive Orlicz spaces $L^{\Phi}(\Omega)$;
\item[(c)] variable Lebesgue spaces $L^{p(\cdot)}(\Omega)$ with $p(\cdot)\in \B_{M}(\R)$.
\end{itemize}
These spaces constitute generalizations of the more familiar $L^{p}(\Omega)$ spaces, and their technical definitions will be postponed to Section \ref{sec: preliminaries}. The hypotheses placed on each of them and their importance has already been clarified in previous works. We redirect the interested reader to \cite{KV24} for a more detailed exposition.

%%%%%%%%%%%%%%%%%%%%%%%%%%%%%%%%%%%%%%%%%%%%%%%%%%%%%%%%%%%%%%%%%%%%%%%%%%%%

Now that the spaces have been revealed, we are left with introducing the operators in question. In order to do that, we are going require a key ingredient: the Fourier transform. For $f\in L^{1}(\R)$, let $Ff$ denote the Fourier transform
\[
(Ff)(\xi) := \int_{- \infty}^{\infty} f(x) e^{i\xi x} \, dx, 
\qquad 
\xi\in \R.
\]
It is of common knowledge that the operator $F$ extends to a bounded linear operator of $L^{2}(\R)$ onto $L^{2}(\R)$, which we also denote by $F$. The inverse of $F$ is given by $(F^{- 1}f)(x) = (2\pi)^{- 1} (Ff)(-x)$ for a.e. $x\in \R$. After this, we can define Fourier multipliers on the space $X(\R)$ (described in (a)-(c)) to be all functions $a\in L^{\infty}(\R)$ for which the densely defined operator
\[
W^{0}(a)f := F^{- 1}(a \cdot Ff),
\]
maps the dense set $L^{2}(\R) \cap X(\R)$ (see, e.g., \cite[Lemma~2.2]{KS21}) into $X(\R)$ and extends to a bounded linear operator on $X(\R)$. The operator $W^{0}(a)$ is called a Fourier convolution operator and the function $a$ will also be referred to as the symbol of $W^{0}(a)$. It can be readily checked that the set of all Fourier multipliers on $X(\R)$, denoted by $\M_{X}$, is a unital normed algebra under pointwise operations and the norm
\[
\| a \|_{\M_{X}} := \| W^{0}(a) \|_{\B(X(\R))}.
\]
In fact, the properties of the underlying function space $X(\R)$ will make it so that $\M_{X}$ is a Banach algebra (see, e.g., \cite[Corollary~2.4]{KS21} and also \cite[Section~2]{KV24}). Before moving on, it is fundamental to make sure that $\M_{X}$ is non-empty. The answer to this predicament will require the use of functions of bounded variation. We say that a function $a: \R \to \C$ is of bounded variation on $\R$ if the (possibly infinite) number
\[
V(a) := \sup \sum_{j = 1}^{n} |a(x_{j}) - a(x_{j - 1})|,
\]
where the supremum is taken over all partitions of $\R$ of the form
\[
- \infty < x_{0} < x_{1} < \dots < x_{n} < + \infty,
\]
is finite. The set $BV(\R)$ of all functions of bounded variation on 
$\mathbb{R}$ with the norm
\[
\| a \|_{BV(\R)} := \| a \|_{L^{\infty}(\R)} + V(a),
\]
is a unital non-separable Banach algebra. According to \cite[Theorem~4.3]{K15} (see also \cite[Section~2]{KV24}), each $a\in BV(\R)$ is a Fourier multiplier on $X(\R)$ with
\[
\| W^{0}(a) \|_{\B(X(\R))} \le c_{X} \| a \|_{BV(\R)},
\]
where $c_{X}$ is a positive constant depending only on $X(\R)$. We will later use this result to define the class of symbols that we will consider. The operator of interest is not too far out of reach, in fact, the only thing we require are two auxiliary yet simple maps: the restriction operator $r_{+}$ from $X(\R)$ into $X(\R_{+})$, and the zero extension operator $\ell_{+}$ from $X(\R_{+})$ into $X(\R)$. Under these conditions, given $a\in \M_{X}$, we define the Wiener-Hopf operator with symbol $a$ by
\[
W(a) := r_{+} W^{0}(a) \ell_{+}.
\]
The only thing left to introduce are the symbols that we are going to study, which will be the continuous ones. For a subset $S$ of a Banach space $\X$, let $\clos_{\X}(S)$ denote the closure of $S$ with respect to the norm of $\X$. Consider the following Banach algebra of Fourier multipliers:
\[
C_{X}(\dot{\R}) := \clos_{\M_{X}} (BV(\R) \cap C(\dot{\R})),
\]
where $C(\dot{\R})$ stands for the $C^{\ast}$-algebra of continuous functions on the one-point compactification $\dot{\R} = \R \cup \{\infty\}$ of the real line. The estimate established in \cite[Theorem~2.3]{KS21} (combined with explanations in \cite[Section~2]{KV24}) yields the inclusion $C_{X}(\dot{\R}) \subseteq C(\dot{\R})$. Because of this, we refer to $C_{X}(\dot{\R})$ as the algebra of continuous Fourier multipliers on the space $X(\R)$.

%%%%%%%%%%%%%%%%%%%%%%%%%%%%%%%%%%%%%%%%%%%%%%%%%%%%%%%%%%%%%%%%%%%%%%%%%%%%

The Fredholm theory of Wiener-Hopf operators was started by Kre\u{\i}n in his fundamental paper \cite{K62}. His ideas were then improved upon in his collaboration with Gohberg \cite{GK60}, and later by Gohberg and Fel'dman \cite{GF74}. These last two works determined the present face of the theory of Wiener-Hopf operators. The next significant step in the study of Wiener-Hopf operators was done by Duduchava in \cite{D79}, who considered discontinuous symbols for the first time, culminating in \cite[Theorem~4.2]{D79}. Later on, Duduchava's result was extended by B\"{o}ttcher and Spitkovsky in \cite{BS93} to the setting of weighted Lebesgue spaces. For the current state of the art of the theory of Wiener-Hopf operators and their weighted analogues, we refer to the monographs \cite{D79}, \cite{BK97}, \cite{BKS02}, \cite[Chap.~9]{BS06} and \cite[Chap.~5]{RSS11}.

%%%%%%%%%%%%%%%%%%%%%%%%%%%%%%%%%%%%%%%%%%%%%%%%%%%%%%%%%%%%%%%%%%%%%%%%%%%%

The aim of this paper is twofold: to shed light on Duduchava's work on Wiener-Hopf operators, specifically with continuous symbols, by providing a detailed proof of \cite[Theorem~3.2]{D79}, and to showcase the limits of his approach by pushing it to other classes of function spaces. The main result says the following.

%%%%%%%%%%%%%%%%%%%%%%%%%%%%%%%%%%%%%%%%%%%%%%%%%%%%%%%%%%%%%%%%%%%%%%%%%%%%

\begin{theorem}\label{thm: main result}
Let $X(\R_{+})$ be one of the three Banach function spaces indicated in (a), (b), (c) above, and $a\in C_{X}(\dot{\R})$. For the Wiener-Hopf operator $W(a)$ to be Fredholm on $X(\R_{+})$ it is necessary and sufficient that the symbol $a$ be elliptic:
\[
\forall \xi\in \dot{\R} 
\quad 
a(\xi) \neq 0.
\]
If the Wiener-Hopf operator $W(a)$ is Fredholm on $X(\R_{+})$, then
\[
\Ind W(a) = - \wind a,
\]
where $\wind a$ is the winding number of $a$ around the origin.
\end{theorem}

%%%%%%%%%%%%%%%%%%%%%%%%%%%%%%%%%%%%%%%%%%%%%%%%%%%%%%%%%%%%%%%%%%%%%%%%%%%%

The paper is organized as follows. In Section \ref{sec: preliminaries}, we gather the definitions of Lorentz, Orlicz and variable Lebesgue spaces. Afterwards, we collect a few results about interpolation which will, in the upcoming sections, be of great importance for the purpose of transferring properties of operators, namely, the boundedness and compactness to each of these function spaces. For each of these properties, there are two results for a simple reason: in contrast with Lorentz and Orlicz spaces, variable Lebesgue spaces are not rearrangement-invariant.

In Section \ref{sec: algebra C_X(R)}, we show that the set of all continuous piecewise linear functions on $\dot{\R}$ is dense in $C_{X}(\dot{\R})$ (with respect to the norm in $\M_{X}$) in each of the aforementioned function spaces by means of interpolation. After this, we establish that the inverse of every elliptic continuous symbol is also a continuous symbol. This is obtained by computing the maximal ideal space of $C_{X}(\dot{\R})$ followed by a direct application of Gelfand's representation theorem.
 
Finally, in Section \ref{sec: final prep}, we provide a self-contained proof of the compactness of the commutator $[aI, S_{\R}]$, where $a\in C(\dot{\R})$ and $S_{\R}$ is the Cauchy singular integral operator on the real line. This result is then applied to establish the compactness of semi-commutators between two Wiener-Hopf operators with continuous symbols. Once this is done, we employ some recent results from \cite{KV24} to derive the Fredholmness of Wiener-Hopf operators generated by an important class of rational symbols. After this, we establish a homotopy between an elliptic continuous piecewise linear function and a suitable rational symbol. Armed with these results, we prove Theorem \ref{thm: main result}.

%%%%%%%%%%%%%%%%%%%%%%%%%%%%%%%%%%%%%%%%%%%%%%%%%%%%%%%%%%%%%%%%%%%%%%%%%%%%
%%%%%%%%%%%%%%%%%%%%%%%%%%%%%%%%%%%%%%%%%%%%%%%%%%%%%%%%%%%%%%%%%%%%%%%%%%%%
%%%%%%%%%%%%%%%%%%%%%%%%%%%%%%%%%%%%%%%%%%%%%%%%%%%%%%%%%%%%%%%%%%%%%%%%%%%%

\section{Preliminaries}\label{sec: preliminaries}
\subsection{Lorentz spaces}

Let $\Omega\in \left\{\R_{+}, \R\right\}$. The distribution function of an a.e. finite Lebesgue measurable function $f: \Omega \to \C$ is given by
\[
\mu_{f}(\lambda) := |\left\{x\in \Omega : |f(x)| > \lambda\right\}|, 
\qquad 
\lambda\geq 0,
\]
where $|E|$ denotes the Lebesgue measure of a measurable subset $E \subseteq \Omega$. The non-increasing rearrangement of $f$ is the function defined by
\[
f^{\ast}(t) := \inf\left\{\lambda\geq 0 : \mu_{f}(\lambda)\le t\right\},
\qquad 
t\geq 0,
\]
where we use the standard convention that $\inf \emptyset = + \infty$.

For $0 < p, q \le \infty$, the Lorentz space $L^{p, q}(\Omega)$ consists of all Lebesgue measurable functions $f: \Omega \to \C$ such that the quantity
\[
\| f \|_{L^{p, q}(\Omega)} 
:= 
\begin{cases}
\bigg(\displaystyle\int_{0}^{\infty} (t^{1/p} f^{\ast \ast}(t))^{q} 
\dfrac{dt}{t}\bigg)^{1/q}, &\ 0 < q < \infty, \\
\displaystyle\sup_{t > 0} t^{1/p} f^{\ast \ast}(t), &\ q = \infty,
\end{cases}
\]
is finite, where
\[
f^{\ast \ast}(t) := \dfrac{1}{t} \int_{0}^{t} f^{\ast}(s) \ ds, 
\qquad 
t > 0.
\]
According to \cite[Chap.~4, Theorem~4.6]{BS88}, if $1 < p < \infty$ and $1\le q\le \infty$, then $(L^{p, q}(\Omega), \| \cdot \|_{L^{p, q}(\Omega)})$ is a rearrangement-invariant Banach function space in the sense of \cite[Chap.~2]{BS88}.

By \cite[Chap.~4, Corollary~4.8]{BS88} (see also \cite[Corollary~8.5.4]{PKJF13}), the space $L^{p, q}(\Omega)$ is separable provided $1 < p < \infty$ and $1\le q < \infty$. As a result of this and \cite[Theorem~8.3.1]{PKJF13}, it follows that $L^{p, q}(\Omega)$ is reflexive if $1 < p, q < \infty$ (see \cite[Corollary~8.5.5]{PKJF13}).

%%%%%%%%%%%%%%%%%%%%%%%%%%%%%%%%%%%%%%%%%%%%%%%%%%%%%%%%%%%%%%%%%%%%%%%%%%%%

\subsection{Orlicz spaces}

Let $\phi: [0, \infty) \to [0, \infty]$ be an increasing left-continuous function satisfying $\phi(0) = 0$ and $\phi(\infty) = \infty$. The function $\Phi$ defined by
\[
\Phi(x) := \int_{0}^{x} \phi(t) \ dt, 
\qquad 
x\geq 0,
\]
is called a Young function. 

Let $\Omega\in \left\{\R_{+}, \R\right\}$. Given a Young function $\Phi$, the Orlicz space $L^{\Phi}(\Omega)$ consists of all Lebesgue measurable functions $f: \Omega \to \C$ such that
\[
\int_{\Omega} \Phi\bigg(\dfrac{|f(x)|}{\lambda}\bigg) \ dx < \infty
\]
for some $\lambda = \lambda(f) > 0$. By \cite[Chap.~4, Theorem~8.9]{BS88}, 
the space $L^{\Phi}(\Omega)$ endowed with the norm
\[
\| f \|_{L^{\Phi}(\Omega)} := \inf\left\{\lambda > 0 : \int_{\Omega} \Phi\bigg(\dfrac{|f(x)|}{\lambda}\bigg) \ dx \le 1\right\}, 
\]
is a rearrangement-invariant Banach function space in the sense of \cite[Chap.~2]{BS88}.

%%%%%%%%%%%%%%%%%%%%%%%%%%%%%%%%%%%%%%%%%%%%%%%%%%%%%%%%%%%%%%%%%%%%%%%%%%%%

\subsection{Variable Lebesgue spaces}

Let $\Omega\in \left\{\R_{+}, \R\right\}$. The set of all measurable a.e. finite functions $p(\cdot): \Omega\to [1, \infty]$ will be denoted by $\mathcal{P}(\Omega)$. The elements of $\mathcal{P}(\Omega)$ are called variable exponents. For $p(\cdot)\in \mathcal{P}(\Omega)$, the variable Lebesgue space $L^{p(\cdot)}(\Omega)$ consists of all Lebesgue measurable functions $f: \Omega \to \C$ such that
\[
\int_{\Omega} \bigg|\dfrac{f(x)}{\lambda}\bigg|^{p(x)} dx < \infty
\]
for some $\lambda = \lambda(f) > 0$. By \cite[Theorem~2.71]{CF13},  
$L^{p(\cdot)}(\Omega)$ is a Banach space with the norm
\[
\| f \|_{L^{p(\cdot)}(\Omega)} := \inf\left\{\lambda > 0 : \int_{\Omega} \bigg|\dfrac{f(x)}{\lambda}\bigg|^{p(x)} dx \le 1\right\}.
\]
Denote
\[
p_{-} := \essinf_{x\in \Omega} p(x), 
\qquad 
p_{+}:= \esssup_{x\in \Omega} p(x).
\]
Obviously, $1\le p_{-} \le p_{+} \le \infty$. It is observed in \cite[Section~2.10.3]{CF13} that $L^{p(\cdot)}(\Omega)$ is a Banach function space in the sense of \cite[Chap.~1]{BS88}. According to \cite[Corollary~2.81]{CF13}, the space $L^{p(\cdot)}(\Omega)$ is reflexive if and only if
\begin{equation}\label{eq: reflexivity in Lp() spaces}
1 < p_{-} \le p_{+} < \infty.
\end{equation}
We shall supplement condition \eqref{eq: reflexivity in Lp() spaces} with the assumption that the Hardy-Littlewood maximal operator $M$ defined by
\[
(Mf)(x) := \sup_{I\ni x} \dfrac{1}{|I|} \int_{I} |f(t)| \ dt,
\]
where the supremum is taken over all bounded intervals $I$ containing the point $x\in \R$, be bounded on $L^{p(\cdot)}(\R)$. The set of all variable exponents satisfying these conditions will be denoted by $\B_{M}(\R)$.

%%%%%%%%%%%%%%%%%%%%%%%%%%%%%%%%%%%%%%%%%%%%%%%%%%%%%%%%%%%%%%%%%%%%%%%%%%%%

\subsection{Interpolation toolkit}

Let $\Omega\in \{\R_{+},\R\}$ and let $X(\Omega)$ be a rearrangement-invariant Banach function space (see \cite[Chap.~2]{BS88}) with the Boyd indices $\alpha_{X}, \beta_{X}$. The class of rearrangement-invariant Banach function spaces is very wide. It includes all Lebesgue spaces $L^{p}(\Omega)$ with $1\le p \le \infty$, all Orlicz spaces $L^\Phi(\R)$, and all Lorentz spaces $L^{p,q}(\Omega)$ with $1< p <\infty$ and $1\le q \le\infty$. We refer to \cite{B69}, \cite[Chap.~3, Section~5]{BS88} and also to \cite[Section~2]{KV24} for the definition and basic properties of Boyd indices. These indices satisfy $0 \le \alpha_{X} \le \beta_{X} \le 1$ and are generalizations of the number $1/p$ for $L^p(\Omega)$ with $1 \le p \le \infty$. The Boyd indices are said to be nontrivial if $0 < \alpha_{X} \le \beta_{X} < 1$. The Boyd indices of the Lorentz space $L^{p,q}(\Omega)$ with $1 < p < \infty$ and $1 \le q \le \infty$ are both equal to $1/p$ (see \cite[Chap.~4, Theorem 4.6]{BS88}). It is also known that an Orlicz space $L^\Phi(\Omega)$ is reflexive if and only if its Boyd indices are nontrivial (see, e.g., \cite[Theorem~3]{KV24}).

%%%%%%%%%%%%%%%%%%%%%%%%%%%%%%%%%%%%%%%%%%%%%%%%%%%%%%%%%%%%%%%%%%%%%%%%%%%%

One key tool that we are going to need throughout the paper is interpolation. Interpolation is a powerful tool that enables us to transfer properties of operators such as boundedness and compactness. We first focus on boundedness and begin by presenting a result which follows from Boyd's interpolation theorem \cite[Theorem~1]{B69}.

%%%%%%%%%%%%%%%%%%%%%%%%%%%%%%%%%%%%%%%%%%%%%%%%%%%%%%%%%%%%%%%%%%%%%%%%%%%%

\begin{theorem}\label{thm: interpolation, boundedness, r.i. BFS}
Let $\Omega\in \left\{\R_{+}, \R\right\}$, $1 \le p < q \le \infty$ and $X(\Omega)$ be a rearrangement-invariant Banach function space whose Boyd indices $\alpha_{X}$ and $\beta_{X}$ satisfy
\begin{equation}\label{eq: interpolation condition on Boyd indices}
\dfrac{1}{q} < \alpha_{X} \le \beta_{X} < \dfrac{1}{p}.
\end{equation}
Under these conditions, there exists a constant $C_{p, q} > 0$ with the following property: if a linear operator $T$ is bounded on $L^{p}(\Omega)$ and $L^{q}(\Omega)$, then it is also bounded on $X(\Omega)$ and
\[
\| T \|_{\B(X(\Omega))} \le C_{p, q} \max\left\{\| T \|_{\B(L^{p}(\Omega))}, \| T \|_{\B(L^{q}(\Omega))}\right\}.
\]
\end{theorem}

%%%%%%%%%%%%%%%%%%%%%%%%%%%%%%%%%%%%%%%%%%%%%%%%%%%%%%%%%%%%%%%%%%%%%%%%%%%%

Unfortunately, variable Lebesgue spaces do not fit the required pre-requisites of Theorem \ref{thm: interpolation, boundedness, r.i. BFS} since they are not rearrangement-invariant. The variable Lebesgue space counterpart is the following.

%%%%%%%%%%%%%%%%%%%%%%%%%%%%%%%%%%%%%%%%%%%%%%%%%%%%%%%%%%%%%%%%%%%%%%%%%%%%

\begin{theorem}[{\cite[Theorem~1.3]{KS23}}]\label{thm: interpolation, boundedness, variable Lp spaces}
Let $p_{j}(\cdot)\in \mathcal{P}(\R)$, $j = 0, 1$. For each $\theta\in (0, 1)$, consider the variable exponent $p_{\theta}(\cdot)\in \mathcal{P}(\R)$ defined by
\[
\dfrac{1}{p_{\theta}(x)} = \dfrac{1 - \theta}{p_{0}(x)} + \dfrac{\theta}{p_{1}(x)}, 
\qquad 
x\in \R.
\]
If a linear mapping
\[
T: L^{p_{0}(\cdot)}(\R) + L^{p_{1}(\cdot)}(\R) \to L^{p_{0}(\cdot)}(\R) + L^{p_{1}(\cdot)}(\R),  
\]
is such that the restriction of $T$ is bounded on $L^{p_{j}(\cdot)}(\R)$ for $j = 0, 1$, then the restriction of $T$ is bounded on $L^{p_{\theta}(\cdot)}(\R)$ and
\[
\| T \|_{\B(L^{p_{\theta}(\cdot)}(\R))} \le 2 \| T \|_{\B(L^{p_{0}(\cdot)}(\R))}^{1 - \theta} \| T \|_{\B(L^{p_{1}(\cdot)}(\R))}^{\theta}.
\]
\end{theorem}

%%%%%%%%%%%%%%%%%%%%%%%%%%%%%%%%%%%%%%%%%%%%%%%%%%%%%%%%%%%%%%%%%%%%%%%%%%%%

As previously stated, the second property of interest to us in interpolating is compactness. Here are the results that we are going to need.

%%%%%%%%%%%%%%%%%%%%%%%%%%%%%%%%%%%%%%%%%%%%%%%%%%%%%%%%%%%%%%%%%%%%%%%%%%%%

\begin{theorem}[{\cite[Corollary~1]{S72}}]\label{thm: interpolation, compactness, r.i. BFS}
Let $\Omega\in \{\R_{+}, \R\}$, $1 < p < q < \infty$ and $T$ be a linear operator that is bounded on $L^{p}(\Omega)$ and $L^{q}(\Omega)$. If $T$ is compact on $L^{p}(\Omega)$ or on $L^{q}(\Omega)$, then it is compact on every rearrangement-invariant Banach function space $X(\Omega)$ whose Boyd indices $\alpha_{X}$ and $\beta_{X}$ satisfy \eqref{eq: interpolation condition on Boyd indices}.
\end{theorem}

%%%%%%%%%%%%%%%%%%%%%%%%%%%%%%%%%%%%%%%%%%%%%%%%%%%%%%%%%%%%%%%%%%%%%%%%%%%%

\begin{theorem}[{\cite[Lemma~6.4]{K15}}]\label{thm: interpolation, compactness, variable Lp spaces}
Let $T\in \B(L^{p(\cdot)}(\R))$ for all $p(\cdot)\in \B_{M}(\R)$. If $T\in \K(L^{r}(\R))$ for some $r\in (1, \infty)$, then $T\in \K(L^{p(\cdot)}(\R))$ for all $p(\cdot)\in \B_{M}(\R)$.
\end{theorem}

%%%%%%%%%%%%%%%%%%%%%%%%%%%%%%%%%%%%%%%%%%%%%%%%%%%%%%%%%%%%%%%%%%%%%%%%%%%%
%%%%%%%%%%%%%%%%%%%%%%%%%%%%%%%%%%%%%%%%%%%%%%%%%%%%%%%%%%%%%%%%%%%%%%%%%%%%
%%%%%%%%%%%%%%%%%%%%%%%%%%%%%%%%%%%%%%%%%%%%%%%%%%%%%%%%%%%%%%%%%%%%%%%%%%%%

\section{Algebra \texorpdfstring{$C_{X}(\dot{\R})$}{TEXT}}\label{sec: algebra C_X(R)}
\subsection{Algebra \texorpdfstring{$C_{X}(\dot{\R})$}{TEXT} as the closure of the set of continuous piecewise linear functions}\label{subsec: piecewise linear}

A function $a\in L^\infty(\R)$ is called piecewise linear if there is a partition of $\R$ with vertices
\[
-\infty < x_{1} < \dots < x_{n} < \infty
\]
and complex constants $c_{k}$, $d_{k}\, (k = 0, \dots, n - 1)$ such that
\[
a(\xi) = c_{0} \chi_{(-\infty, x_{1}]}(\xi) + \sum_{k = 1}^{n - 1} (c_{k} + d_{k} \xi) \chi_{(x_{k}, x_{k+1}]}(\xi) + d_{0} \chi_{(x_{n}, \infty)}(\xi).
\]
The set of all piecewise linear functions will be denoted by $P\ell$.

%%%%%%%%%%%%%%%%%%%%%%%%%%%%%%%%%%%%%%%%%%%%%%%%%%%%%%%%%%%%%%%%%%%%%%%%%%%%

A routine application of the definition of bounded variation yields the following.

%%%%%%%%%%%%%%%%%%%%%%%%%%%%%%%%%%%%%%%%%%%%%%%%%%%%%%%%%%%%%%%%%%%%%%%%%%%%

\begin{lemma}\label{lem: Pell subset of BV}
$P\ell$ is a subset of $BV(\R)$.
\end{lemma}

%%%%%%%%%%%%%%%%%%%%%%%%%%%%%%%%%%%%%%%%%%%%%%%%%%%%%%%%%%%%%%%%%%%%%%%%%%%%

%\subsection{The case of separable rearrangement-invariant spaces.}

For each $1 < p < \infty$, we shall denote $\M_{p} := \M_{L^{p}}$.

%%%%%%%%%%%%%%%%%%%%%%%%%%%%%%%%%%%%%%%%%%%%%%%%%%%%%%%%%%%%%%%%%%%%%%%%%%%%

\begin{theorem}\label{thm: C_X in terms of piecewise linear - RI BFS version}
For every separable rearrangement-invariant Banach function space $X(\R)$ with nontrivial Boyd indices,
\[
C_{X}(\dot{\R}) = \clos_{\M_{X}} (P\ell \cap C(\dot{\R})).
\]
\end{theorem}

\begin{proof}
The proof is similar to that of \cite[Theorem~3.3]{FKV21}. By Lemma \ref{lem: Pell subset of BV}, we have
\[
P\ell \cap C(\dot{\R}) \subseteq BV(\R) \cap C(\dot{\R}) \subseteq C_{X}(\dot{\R}).
\]
Bearing this in mind, let $a\in C_{X}(\dot{\R})$ and fix $\varepsilon > 0$. By definition of $ C_{X}(\dot{\R})$, there is some $b\in BV(\R)\cap C(\dot{\R})$ such that
\[
\| a - b \|_{\M_{X}} < \dfrac{\varepsilon}{2}.
\]
According to \cite[Lemma~2.10]{D79}, there is a sequence $(c_{n})_{n\in \N}$ in $P\ell \cap C(\dot{\R})$ such that
\[
\sup_{n\in \N} V(c_{n}) \le V(b), \qquad \lim_{n\to \infty} \| b - c_{n} \|_{L^{\infty}(\R)} = 0.
\]
With this in mind, consider
\[  
p := \bigg(\dfrac{\beta_{X} + 1}{2}\bigg)^{- 1},
\qquad
q := \bigg(\dfrac{\alpha_{X}}{2}\bigg)^{- 1}.
\]
The assumption that $\alpha_{X}, \beta_{X}\in (0, 1)$ gives $1 < p, q < \infty$. Following the procedure presented in \cite[Remark~2.12]{D79} (see also \cite[Theorem~3.2]{FKV21}), we can find some $n_{0}\in \N$ sufficiently large so that
\[
\| b - c_{n_{0}} \|_{\M_{p}} < \dfrac{\varepsilon}{2 C_{p, q}}, 
\qquad 
\| b - c_{n_{0}} \|_{\M_{q}} < \dfrac{\varepsilon}{2 C_{p, q}},
\]
where $C_{p, q}$ is the constant from Theorem \ref{thm: interpolation, boundedness, r.i. BFS}. Under these conditions, since by construction we have
\[
\dfrac{1}{q} < \alpha_{X} \le \beta_{X} < \dfrac{1}{p},
\]
a direct application of Theorem \ref{thm: interpolation, boundedness, r.i. BFS} yields
\begin{align*}
\| b - c_{n_{0}} \|_{\M_{X}} &= \| W^{0}(b - c_{n_{0}}) \|_{\B(X(\R))} \\ 
&\le C_{p, q} \max\left\{\| W^{0}(b - c_{n_{0}}) \|_{\B(L^{p}(\R))}, \| W^{0}(b - c_{n_{0}}) \|_{\B(L^{q}(\R))}\right\} \\
&= C_{p, q} \max\left\{\| b - c_{n_{0}} \|_{\M_{p}}, \| b - c_{n_{0}} \|_{\M_{q}}\right\} \\
&< \dfrac{\varepsilon}{2}.
\end{align*}
Combining these remarks, we get
\[
\| a - c_{n_{0}} \|_{\M_{X}} \le \| a - b \|_{\M_{X}} + \| b - c_{n_{0}} \|_{\M_{X}} < \varepsilon.
\]
Since $\varepsilon > 0$ was considered arbitrary, the result follows.
\end{proof}

%%%%%%%%%%%%%%%%%%%%%%%%%%%%%%%%%%%%%%%%%%%%%%%%%%%%%%%%%%%%%%%%%%%%%%%%%%%%

%\subsection{The case of variable Lebesgue spaces.}

Given $p(\cdot)\in \B_{M}(\R)$, we shall use the following abbreviations:
\[
\M_{p(\cdot)} := \M_{L^{p(\cdot)}}, 
\qquad 
C_{L^{p(\cdot)}}(\dot{\R}) := C_{p(\cdot)}(\dot{\R}). 
\]

%%%%%%%%%%%%%%%%%%%%%%%%%%%%%%%%%%%%%%%%%%%%%%%%%%%%%%%%%%%%%%%%%%%%%%%%%%%%

\begin{theorem}\label{thm: C_p() in terms of piecewise linear}
For all $p(\cdot)\in \B_{M}(\R)$,
\[
C_{p(\cdot)}(\dot{\R}) = \clos_{\M_{p(\cdot)}} (P\ell \cap C(\dot{\R})).
\]
\end{theorem}

\begin{proof}
The proof follows the same argument as presented in \cite[Theorem~4.2]{K20}. By Lemma \ref{lem: Pell subset of BV}, it follows immediately that
\[
P\ell \cap C(\dot{\R}) \subseteq C_{p(\cdot)}(\dot{\R}).
\]
Bearing this in mind, let $a\in C_{p(\cdot)}(\dot{\R})$ and fix $\varepsilon > 0$. By definition of $C_{p(\cdot)}(\dot{\R})$, there exists some $b\in BV(\R) \cap C(\dot{\R})$ such that
\begin{equation}\label{eq: estimate of a - b in M_p()}
\| a - b \|_{\M_{p(\cdot)}} < \dfrac{\varepsilon}{2}.
\end{equation}
In light of \cite[Lemma~2.10]{D79}, there is a sequence $(c_{n})_{n\in \N}$ in $P\ell \cap C(\dot{\R})$ such that for all $n\in \N$,
\[
V(c_{n}) \le V(b), 
\qquad 
\lim_{n\to \infty} \| b - c_{n} \|_{L^{\infty}(\R)} = 0.
\]
Then for large enough $N\in \N$ we have
\[
\| b - c_{n} \|_{L^{\infty}(\R)} < \| b \|_{BV(\R)}, 
\qquad 
n > N.
\]
From here it follows that for all $n > N$,
\[
\| c_{n} \|_{L^\infty(\R)} - \| b \|_{L^\infty(\R)} \le \big| \| b \|_{L^\infty(\R)} - \| c_{n} \|_{L^\infty(\R)} \big| \le \| b - c_{n} \|_{L^\infty(\R)} < \| b \|_{BV(\R)},
\]
and, consequently,
\[
\| c_{n} \|_{L^{\infty}(\R)} < \| b \|_{BV(\R)} + \| b \|_{L^{\infty}(\R)}.
\]
Adding the term $V(c_{n})$ on both sides of the inequality and using the definition of $\| \cdot \|_{BV(\R)}$ we arrive at
\[
\| c_{n} \|_{BV(\R)} < \| b \|_{BV(\R)} + \| b \|_{L^\infty(\R)} + V(c_{n}),
\qquad 
n > N.
\]
Finally, applying the estimate $V(c_{n}) \le V(b)$ and again the definition of $\| \cdot \|_{BV(\R)}$ we get
\begin{equation}\label{eq: preliminar estimate of norm in BV of c_n in terms of b}
\| c_{n} \|_{BV(\R)} < 2 \| b \|_{BV(\R)}, 
\qquad 
n > N.
\end{equation}
We will now set up for an interpolation argument: since $p(\cdot)\in \B_{M}(\R)$, by \cite[Theorem~1.3]{DKS22}, there exists a number $\Theta_{p(\cdot), 2}\in (0, 1)$ such that for all $\theta\in (0, \Theta_{p(\cdot), 2}]$ the variable exponent $p_{\theta}(\cdot)$ defined by
\begin{equation}\label{eq: interpolation in variable lebesgue space}
\dfrac{1}{p(x)} = \dfrac{\theta}{2} + \dfrac{1 - \theta}{p_{\theta}(x)},
\end{equation}
belongs to $\B_{M}(\R)$. Set $\theta := \Theta_{p(\cdot), 2}$. It follows from the Stechkin-type inequality \cite[Theorem~4.3]{K15} and inequality \eqref{eq: preliminar estimate of norm in BV of c_n in terms of b} that for all $n > N$,
\begin{equation}\label{eq: estimate of norm in M_ptheta() of c_n - b in terms of b}
\| c_{n} - b \|_{\M_{p_{\theta}(\cdot)}} \le \| c_{n} \|_{\M_{p_{\theta}(\cdot)}} + \| b \|_{\M_{p_{\theta}(\cdot)}} \le 3 \| b \|_{\M_{p_{\theta}(\cdot)}} \le 3c_{L^{p_{\theta}(\cdot)}} \| b \|_{BV(\R)}.
\end{equation}
Identity \eqref{eq: interpolation in variable lebesgue space} combined with Theorem \ref{thm: interpolation, boundedness, variable Lp spaces} and inequality \eqref{eq: estimate of norm in M_ptheta() of c_n - b in terms of b} yield that for every $n > N$,
\begin{align}
\| c_{n} - b \|_{\M_{p(\cdot)}} = \| W^{0}(c_{n} - b) \|_{\B(L^{p(\cdot)}(\R))} &\le 2 \| W^{0}(c_{n} - b) \|_{\B(L^{2}(\R))}^{\theta} \| W^{0}(c_{n} - b) \|_{\B(L^{p_{\theta}(\cdot)}(\R))}^{1 - \theta} \notag \\ 
&= 2 \| c_{n} - b \|_{L^{\infty}(\R)}^{\theta} \| c_{n} - b \|_{\M_{p_{\theta}(\cdot)}}^{1 - \theta} \notag \\
&\le 2 (3c_{L^{p_{\theta}(\cdot)}})^{1 - \theta} \| b \|_{BV(\R)}^{1 - \theta} \| c_{n} - b \|_{L^{\infty}(\R)}^{\theta}. \label{eq: estimate of norm in M_p() of c_n - b in terms of Linfty norm of c_n - b}
\end{align}
It follows from $\| b - c_{n} \|_{L^{\infty}(\R)}\to 0$ and inequality \eqref{eq: estimate of norm in M_p() of c_n - b in terms of Linfty norm of c_n - b} that there exists $n_{0} > N$ such that
\begin{equation}\label{eq: estimate of c_n0 - b in M_p()}
\| c_{n_{0}} - b \|_{\M_{p(\cdot)}} < \dfrac{\varepsilon}{2}.
\end{equation}
Combining inequalities \eqref{eq: estimate of a - b in M_p()} and \eqref{eq: estimate of c_n0 - b in M_p()} we see that there exists a continuous piecewise linear function $c_{n_0}$ such that 
\[
\| a - c_{n_{0}} \|_{\M_{p(\cdot)}} < \varepsilon.
\]
Since $\varepsilon > 0$ was considered arbitrary, the result follows.
\end{proof}

%%%%%%%%%%%%%%%%%%%%%%%%%%%%%%%%%%%%%%%%%%%%%%%%%%%%%%%%%%%%%%%%%%%%%%%%%%%%

\subsection{\texorpdfstring{$C_{X}(\dot{\R})$}{TEXT} is inverse-closed in the algebra \texorpdfstring{$L^{\infty}(\R)$}{TEXT}}\label{subsec: inverse closedness}

It is well-known (see, e.g., \cite[Pages 811-812]{GGK93}) that the maximal ideal space of $C(\dot{\R})$, denoted by $M(C(\dot{\R}))$, coincides with the set 
\[
\left\{\delta_{x} : x\in \dot{\R}\right\},
\]
where, for each $x\in \dot{\R}$, $\delta_{x} : C(\dot{\R}) \to \C$ corresponds to the evaluation functional defined by $\delta_{x}(f) := f(x)$, $f\in C(\dot{\R})$. It is easily seen that the procedure that took place in achieving this can be repeated to the subalgebra $BV(\R) \cap C(\dot{\R})$ to get that
\[
M(BV(\R) \cap C(\dot{\R})) = \left\{\restr{\delta_{x}}{BV(\R) \, \cap \, C(\dot{\R})} : x\in \dot{\R}\right\}.
\]
As an immediate consequence of this, we arrive at the following conclusion.

%%%%%%%%%%%%%%%%%%%%%%%%%%%%%%%%%%%%%%%%%%%%%%%%%%%%%%%%%%%%%%%%%%%%%%%%%%%%

\begin{lemma}\label{lem: extension of multiplicative linear functionals}
For each functional $\varphi\in M(BV(\R) \cap C(\dot{\R}))$, there is a unique extension $\Phi\in M(C(\dot{\R}))$.
\end{lemma}

%%%%%%%%%%%%%%%%%%%%%%%%%%%%%%%%%%%%%%%%%%%%%%%%%%%%%%%%%%%%%%%%%%%%%%%%%%%%

Once armed with these results, we are ready to compute $M(C_{X}(\dot{\R}))$.

%%%%%%%%%%%%%%%%%%%%%%%%%%%%%%%%%%%%%%%%%%%%%%%%%%%%%%%%%%%%%%%%%%%%%%%%%%%%

\begin{theorem}\label{thm: maximal ideal space of C_X}
If $X(\R)$ is a separable rearrangement-invariant Banach function space with nontrivial Boyd indices, or a variable Lebesgue space with $p(\cdot)\in \mathcal{B}_{M}(\R)$, then the maximal ideal space of $C_{X}(\dot{\R})$ is homeomorphic to $\dot{\R}$.
\end{theorem}

\begin{proof}
Let us first start by pointing out the obvious:
\[
\left\{\restr{\delta_{x}}{C_{X}(\dot{\R})} : x\in \dot{\R}\right\} \subseteq M(C_{X}(\dot{\R})). 
\]
We claim that the converse also holds. To see this, note that we have the following chain of homomorphic embeddings:
\[
BV(\R) \cap C(\dot{\R}) \subseteq C_{X}(\dot{\R}) \subseteq C(\dot{\R}).
\]
By definition of the algebra $C_{X}(\dot{\R})$, the algebra $BV(\R) \cap C(\dot{\R})$ is dense in $C_{X}(\dot{\R})$ with respect to the norm of $\M_{X}$. On account of these remarks and Lemma \ref{lem: extension of multiplicative linear functionals}, we see that the commutative Banach algebras
\[
\mathfrak{A} := BV(\R) \cap C(\dot{\R}), \qquad \mathfrak{B} := C_{X}(\dot{\R}), \qquad \mathfrak{C} := C(\dot{\R}),
\]
satisfy the conditions of \cite[Theorem~3.5]{FKK21}. As a consequence of this, we get that every functional $\varphi\in M(C_{X}(\dot{\R}))$ admits a unique extension $\Phi\in M(C(\dot{\R}))$. Taking into account that
\[
M(C(\dot{\R})) = \left\{\delta_{x} : x\in \dot{\R}\right\},
\]
we find that there is some $x\in \dot{\R}$ satisfying $\Phi = \delta_{x}$, and thus
\[
\varphi = \restr{\delta_{x}}{C_{X}(\dot{\R})}.
\]
This proves our claim. The only thing left to do is to verify that the mapping
\[
\varrho: \dot{\R} \to M(C_{X}(\dot{\R})), 
\qquad 
\varrho(x) := \restr{\delta_{x}}{C_{X}(\dot{\R})},
\]
is a homeomorphism. The bijectivity has already been proved above. Regarding the continuity, fix $x\in \dot{\R}$ and let $(x_{n})_{n\in \N}$ be a sequence of elements in $\dot{\R}$ such that
\[
\lim_{n\to \infty} x_{n} = x.
\]
The continuity of the functions in $C_{X}(\dot{\R})$ implies that for all $a\in C_{X}(\dot{\R})$,
\[
\lim_{n\to \infty} a(x_{n}) = a(x),
\]
or, equivalently,
\[
\lim_{n\to \infty} \restr{\delta_{x_{n}}}{C_{X}(\dot{\R})}(a) = \restr{\delta_{x}}{C_{X}(\dot{\R})}(a),
\]
Since this holds for any $a\in C_{X}(\dot{\R})$, we get that the sequence of functionals
\[
(\restr{\delta_{x_{n}}}{C_{X}(\dot{\R})})_{n\in \N}, 
\]
weak$^{\ast}$-converges to $\restr{\delta_{x}}{C_{X}(\dot{\R})}$ as intended. Therefore $\varrho$ is continuous. Finally, taking into account that $\dot{\R}$ and $M(C_{X}(\dot{\R}))$ are compact Hausdorff spaces (see, e.g., \cite[Theorem~29.1]{M00} and \cite[Theorem~2.1.4]{RSS11}, respectively), a direct application of \cite[Theorem~26.6]{M00} gives the result.
\end{proof}

%%%%%%%%%%%%%%%%%%%%%%%%%%%%%%%%%%%%%%%%%%%%%%%%%%%%%%%%%%%%%%%%%%%%%%%%%%%%

\begin{theorem}\label{thm: inverse closedness of C_X}
Let $X(\R)$ be a separable rearrangement-invariant Banach function space with nontrivial Boyd indices, or a variable Lebesgue space with $p(\cdot)\in \mathcal{B}_{M}(\R)$ and $a\in C_{X}(\dot{\R})$. If the symbol $a$ is elliptic, then $a^{- 1}\in C_{X}(\dot{\R})$.
\end{theorem}

\begin{proof}
Since Theorem \ref{thm: maximal ideal space of C_X} is available, we identify $M(C_{X}(\dot{\R}))$ as $\dot{\R}$ which in turn implies that the Gelfand transform acting on $C_{X}(\dot{\R})$ will simply be the inclusion map onto $C(\dot{\R})$. Faced with this, we are left with appealing to Gelfand's representation theorem (see, e.g., \cite[Theorem~2.1.3, iii)]{RSS11}).
\end{proof}

%%%%%%%%%%%%%%%%%%%%%%%%%%%%%%%%%%%%%%%%%%%%%%%%%%%%%%%%%%%%%%%%%%%%%%%%%%%%
%%%%%%%%%%%%%%%%%%%%%%%%%%%%%%%%%%%%%%%%%%%%%%%%%%%%%%%%%%%%%%%%%%%%%%%%%%%%
%%%%%%%%%%%%%%%%%%%%%%%%%%%%%%%%%%%%%%%%%%%%%%%%%%%%%%%%%%%%%%%%%%%%%%%%%%%%

\section{Proof of the main result}\label{sec: final prep}
\subsection{Compactness of semi-commutators}\label{subsec: compactness of semicommutators}

Let $\Gamma$ be the unit circle $\T$ or the real line $\R$. The Cauchy singular integral operator $S_{\Gamma}$ over the curve $\Gamma$ is defined as
\[
(S_{\Gamma} \varphi)(x) 
:= 
\dfrac{1}{\pi i} \int_{\Gamma} \dfrac{\varphi(t)}{t - x} \ dt,
\qquad
\varphi\in C^{\infty}_{c}(\Gamma),
\]
where the integral is understood in the sense of Cauchy's principal value.

%%%%%%%%%%%%%%%%%%%%%%%%%%%%%%%%%%%%%%%%%%%%%%%%%%%%%%%%%%%%%%%%%%%%%%%%%%%%

The following result is formulated without proof in \cite[Lemma~1.21]{D79}. 
For the sake of completeness, we will give its proof.

%%%%%%%%%%%%%%%%%%%%%%%%%%%%%%%%%%%%%%%%%%%%%%%%%%%%%%%%%%%%%%%%%%%%%%%%%%%%

\begin{theorem}\label{thm: compactness of [aI, S] on Lp}
Let $1 < p < \infty$. For all $a\in C(\dot{\R})$, the commutator between the multiplication operator $aI$ and the Cauchy singular integral operator $S_{\R}$:
\[
[aI, S_{\R}] := aS_{\R} - S_{\R} aI,
\]
is compact on $L^{p}(\R)$.
\end{theorem}

\begin{proof}
Fix $a\in C(\dot{\R})$. Consider the isometric isomorphism between 
$C(\dot{\R})$ and $C(\T)$:
\[
(B_{0}a)(t) 
:= 
\begin{cases}
a\bigg(i \dfrac{1 + t}{1 - t}\bigg), &\text{if} \ t\in \T\setminus \{1\}, \\
a(\infty), &\text{if} \ t = 1,
\end{cases}
\]
alongside the operator
\[
(Bf)(x) := \dfrac{2^{1 - 1/p}}{x + i} f\bigg(\dfrac{x - i}{x + i}\bigg), 
\qquad 
x\in \R.
\]
By \cite[Section~16.1]{BKS02} and \cite[Chap.~1, Theorem~5.1]{GK92}, $B$ defines an isometric isomorphism from the weighted Lebesgue space $L^{p}(\T, w_{p})$,
\[
f\in L^{p}(\T, w_{p}) \ \Leftrightarrow \ fw_{p} \in L^{p}(\T),
\]
where $w_{p}(t) := |t - 1|^{1 - \frac{2}{p}}$, into $L^{p}(\R)$ with inverse
\[
(B^{- 1}g)(t) := i\dfrac{2^{1/p}}{1 - t} g\bigg(i\dfrac{1 + t}{1 - t}\bigg), 
\qquad 
t\in \T\setminus \{1\}.
\]
Furthermore, we have the identity
\[
S_{\R} = B S_{\T} B^{- 1}. 
\]
Combining this with the relation
\[
B^{- 1} aI_{L^{p}(\R)} B = (B_{0}a)I_{L^{p}(\T, w_{p})},
\]
we find that
\[
B^{- 1} [aI_{L^{p}(\R)}, S_{\R}] B = [(B_{0}a)I_{L^{p}(\T, w_{p})}, S_{\T}].
\]
Faced with this, we are left with appealing to \cite[Chap.~1, Theorem~4.3]{GK92} to arrive at the desired result.
\end{proof}

%%%%%%%%%%%%%%%%%%%%%%%%%%%%%%%%%%%%%%%%%%%%%%%%%%%%%%%%%%%%%%%%%%%%%%%%%%%%

\begin{theorem}\label{thm: compactness of semicommutators between WH operators with at least one continuous BV function}
Let $X(\R_{+})$ be one of the three Banach function spaces specified in Theorem \ref{thm: main result} and $a, b\in BV(\R)$. If $a\in C(\dot{\R})$ or $b\in C(\dot{\R})$, then the semi-commutators between the Wiener-Hopf operators $W(a)$ and $W(b)$,
\[
W(a)W(b) - W(ab), 
\qquad 
W(b)W(a) - W(ba),
\]
are compact on $X(\R_{+})$.
\end{theorem}

\begin{proof}
The proof follows the argument made in \cite[Proposition~2.19]{D79}. Fix 
$a, b\in BV(\R)$ and suppose, without loss of generality, that $a\in C(\dot{\R})$. Consider the auxiliary operator
\[
A := \ell_{+} (W(a)W(b) - W(ab)) r_{+}.
\]
It was shown in \cite[Proposition~2.19]{D79} that the operator $A$ on $L^{2}(\R)$ is equal to
\[
\dfrac{1}{4} F^{- 1} [aI, S_{\R}] [bI, S_{\R}] P_{+} F.
\]
The above remarks combined with Theorem \ref{thm: compactness of [aI, S] on Lp} give that $A$ is compact on $L^{2}(\R)$.

If $X = L^{p(\cdot)}$, then since the operator $A$ is compact on $L^{2}(\R)$ and is bounded on $L^{r(\cdot)}(\R)$ for all $r(\cdot)\in \B_{M}(\R)$, Theorem \ref{thm: interpolation, compactness, variable Lp spaces} yields that $A$ is compact on $L^{p(\cdot)}(\R)$. Suppose now that $X = L^{\Phi}$. As an immediate consequence of what we have just proved, the operator $A$ is compact on $L^{r}(\R)$ for all $1 < r < \infty$. Since $L^{\Phi}(\R)$ is reflexive, its Boyd indices are nontrivial (see, e.g., \cite[Theorem~3]{KV24}). Take $p, q\in (1, \infty)$ such that
\[
\dfrac{1}{q} < \alpha_{L^{\Phi}} \le \beta_{L^{\Phi}} < \dfrac{1}{p}.
\]
Under these conditions, since the operator $A$ is compact on $L^{p}(\R)$ and $L^{q}(\R)$, a direct application of Theorem \ref{thm: interpolation, compactness, r.i. BFS} shows that the operator $A$ is compact on $L^{\Phi}(\R)$. In the case when $X = L^{p, q}$, by \cite[Chap.~4, Theorem~4.6]{BS88} one has
\[
\alpha_{L^{p, q}} = \beta_{L^{p, q}} = \dfrac{1}{p}.
\]
Faced with this, one can choose $r := (p + 1)/2$ and $s := p + 1$ to get
\[
0 < \dfrac{1}{s} < \alpha_{L^{p, q}} = \beta_{L^{p, q}} < \dfrac{1}{r} < 1,
\]
and repeat the interpolation argument made for Orlicz spaces.

In either case, we find that the operator $A$ is compact on $X(\R)$. From here we derive that $W(a)W(b) - W(ab)$ is compact on $X(\R_{+})$: indeed, let $(f_{n})_{n\in \N}$ be a bounded sequence in $X(\R_{+})$. The boundedness of $\ell_{+}$ implies that $(\ell_{+} f_{n})_{n\in \N}$ is bounded sequence in $X(\R)$. Employing the compactness of $A$, there exists a subsequence $(\ell_{+} f_{n_{k}})_{k\in \N}$ such that $(A \ell_{+} f_{n_{k}})_{k\in \N}$ converges in $X(\R)$. Consequently, $(r_{+} A \ell_{+} f_{n_{k}})_{k\in \N}$ converges in $X(\R_{+})$. We are left with observing that since $r_{+} \ell_{+} = I_{X(\R_{+})}$, we have
\[
r_{+} A \ell_{+} = W(a)W(b) - W(ab).
\]
\end{proof}

%%%%%%%%%%%%%%%%%%%%%%%%%%%%%%%%%%%%%%%%%%%%%%%%%%%%%%%%%%%%%%%%%%%%%%%%%%%%

Now that Theorem \ref{thm: compactness of semicommutators between WH operators with at least one continuous BV function} is available, a standard density argument yields the following result:

%%%%%%%%%%%%%%%%%%%%%%%%%%%%%%%%%%%%%%%%%%%%%%%%%%%%%%%%%%%%%%%%%%%%%%%%%%%%

\begin{theorem}\label{thm: compactness of semi-commutators between WH operators with continuous symbols}
Let $X(\R_{+})$ be one of the three Banach function spaces specified in Theorem \ref{thm: main result}. For all $a, b\in C_{X}(\dot{\R})$, the semi-commutators between the Wiener-Hopf operators $W(a)$ and $W(b)$,
\[
W(a)W(b) - W(ab), 
\qquad 
W(b)W(a) - W(ba),
\]
are compact on $X(\R_{+})$.
\end{theorem}

\begin{comment}
\begin{proof}
Let $a, b\in C_{X}(\dot{\R})$. According to the definition of the algebra $ C_{X}(\dot{\R})$, there exist sequences $(a_{n})_{n\in \N}$ and $(b_{m})_{m\in \N}$ in $BV(\R) \cap C(\dot{\R})$ such that
\[
\| a - a_{n} \|_{\M_{X}} \to 0, \qquad \| b - b_{m} \|_{\M_{X}} \to 0, 
\]
as $n, m\to \infty$, respectively. Consequently, $ab\in C_{X}(\dot{\R})$ and
\[
\lim_{n, m\to \infty }\| ab - a_{n}b_{m} \|_{\M_{X}} = 0.
\]
Consider the following double sequence of operators
\[
W(a_{n})W(b_{m}) - W(a_{n}b_{m}), \qquad n, m\in \N.
\]
It follows from Theorem \ref{thm: compactness of semicommutators between WH operators with at least one continuous BV function} that
\[
(W(a_{n})W(b_{m}) - W(a_{n}b_{m}))_{n, m\in \N} \subset \K(X(\R_{+})).
\]
With this in mind, straightforward computations show that
\begin{align*}
\| W(a)W(b) - W(ab) &- W(a_{n})W(b_{m}) + W(a_{n}b_{m}) \|_{\B(X(\R_{+}))} \\
&\le \| b \|_{\M_{X}} \| a - a_{n} \|_{\M_{X}} + \sup_{k\in \N} \| a_{k} \|_{\M_{X}} \| b - b_{m} \|_{\M_{X}} + \| ab - a_{n}b_{m} \|_{\M_{X}},
\end{align*}
for all $n, m\in \N$. As a result of this, the above remarks allied with the squeeze theorem yield that
\[
\lim_{n, m\to \infty} W(a_{n})W(b_{m}) - W(a_{n}b_{m}) = W(a)W(b) - W(ab), 
\]
in $\B(X(\R_{+}))$, which in turn reveals that $W(a)W(b) - W(ab)$ is the uniform limit of compact operators on $X(\R_{+})$ and hence also compact.
\end{proof}
\end{comment}

%%%%%%%%%%%%%%%%%%%%%%%%%%%%%%%%%%%%%%%%%%%%%%%%%%%%%%%%%%%%%%%%%%%%%%%%%%%%

\subsection{Fredholmness of \texorpdfstring{$W(r_{n}) \ (n\in \Z)$}{TEXT}}

For each $n\in \Z$, consider the rational function
\[
r_{n}(\xi) := \bigg(\dfrac{\xi - i}{\xi + i}\bigg)^{n}, 
\qquad 
\xi\in \R.
\]
A routine verification shows that $V(r_{n}) = 2\pi|n|$, and so $r_{n}\in C_{X}(\dot{\R})$.

%%%%%%%%%%%%%%%%%%%%%%%%%%%%%%%%%%%%%%%%%%%%%%%%%%%%%%%%%%%%%%%%%%%%%%%%%%%%

\begin{lemma}\label{lem: kernel of W(r_-n)}
Let $X(\R_{+})$ be one of the three Banach function spaces specified in Theorem \ref{thm: main result}. For all $n\in \N$,
\[
\Ker W(r_{- n}) = \Lspan \left\{\psi_{k} : k = 0, 1, \dots, n - 1\right\},
\]
where $\psi_{0}(x) := \sqrt{2} e^{- x}$, $x\in \R_{+}$, and $\psi_{k} := W(r_{k})\psi_{0}$, $k\in \N$.
\end{lemma}

\begin{proof}
Fix $n\in \N$ and consider the auxiliary operator
\[
P_{n} := I - W(r_{n})W(r_{- n}).
\]
It was shown in \cite[Lemma~4]{KV24} that
\[
\Image P_{n} = \Lspan \left\{\psi_{k} : k = 0, 1, \dots, n - 1\right\}.
\]
We claim that $\Ker W(r_{- n}) = \Image P_{n}$. Let us start by observing that $P_{n}^{2} = P_{n}$ and $W(r_{- n}) P_{n} = 0$; indeed,
\begin{align*}
P_{n}^{2} &= P_{n} - W(r_{n})W(r_{- n}) + W(r_{n})W(r_{- n})W(r_{n})W(r_{- n}) \\
&= P_{n} - W(r_{n})W(r_{- n}) + W(r_{n})W(r_{- n}r_{n})W(r_{- n}) \\
&= P_{n},
\end{align*}
and
\[
W(r_{- n}) P_{n} = W(r_{- n}) - W(r_{- n})W(r_{n})W(r_{- n}) = W(r_{- n}) - W(r_{- n}r_{n})W(r_{- n}) = 0,
\]
where in both calculations we employed the identity 
\begin{equation}\label{eq: right invertibility of w(r_-n)}
W(r_{- n})W(r_{n}) = W(r_{- n}r_{n}) = I, 
\end{equation}
which is justified by the use of \cite[Theorem~4]{KV24}. Bearing this in mind, let $f\in X(\R_{+})$. If $f\in \Ker W(r_{- n})$, then $W(r_{- n})f = 0$ which in turn implies that $P_{n}f = f$ and hence $f\in \Image P_{n}$. Conversely, if $f\in \Image P_{n}$, then $P_{n}f = f$ due to $P_{n}$ being idempotent. Applying the Wiener-Hopf operator $W(r_{- n})$ to both sides of this equation, we find that
\[
W(r_{- n})f = W(r_{- n}) P_{n}f = 0, 
\]
and hence $f\in \Ker W(r_{- n})$. Therefore $\Ker W(r_{- n}) = \Image P_{n}$.
\end{proof}

%%%%%%%%%%%%%%%%%%%%%%%%%%%%%%%%%%%%%%%%%%%%%%%%%%%%%%%%%%%%%%%%%%%%%%%%%%%%

\begin{lemma}\label{lem: adjoint of W(a)}
Let $X(\R_{+})$ be one of the three Banach function spaces specified in Theorem \ref{thm: main result}. For all $a\in \M_{X}$, we have $\overline{a}\in \M_{X'}$ and $(W(a))^{\ast} = W(\overline{a})$.
\end{lemma}

\begin{proof}
Fix $a\in \M_{X}$. Since $W(a) := r_{+}W^{0}(a)\ell_{+}$, it suffices to compute the adjoint of each operator separately. It can be readily checked that $(r_{+})^{\ast} = \ell_{+}$ and $(\ell_{+})^{\ast} = r_{+}$. We claim that $(W^{0}(a))^{\ast} = W^{0}(\overline{a})$. To see this, note that for all $f, g\in C^{\infty}_{c}(\R)$, standard properties of the Fourier transform (see, e.g., \cite[Section~11C]{A20}) enable us to write
\[
\int_{\R} (W^{0}(a)f)(x) \overline{g(x)} \ dx 
= 
\int_{\R} f(x) \overline{(W^{0}(\overline{a})g)(x)} \ dx,  
\]
from which the claim follows and, consequently, $\overline{a}\in \M_{X'}$. As a result of this, we get $(W(a))^{\ast} = (\ell_{+})^{\ast} (W^{0}(a))^{\ast} (r_{+})^{\ast} = W(\overline{a})$.
\end{proof}

%%%%%%%%%%%%%%%%%%%%%%%%%%%%%%%%%%%%%%%%%%%%%%%%%%%%%%%%%%%%%%%%%%%%%%%%%%%%

\begin{theorem}\label{thm: Fredholmness of W(r_n)}
Let $X(\R_{+})$ be one of the three Banach function spaces specified in Theorem \ref{thm: main result}. For all $n\in \Z$, the Wiener-Hopf operator $W(r_{n})$ is Fredholm on $X(\R_{+})$ and its Fredholm index is equal to $- n$.
\end{theorem}

\begin{proof}
Fix $n\in \N_{0}$. If $n = 0$, then $W(r_{0}) = I$ which is invertible and hence Fredholm on $X(\R_{+})$ with index $0$. Suppose now that $n\geq 1$. It follows immediately from Equation \eqref{eq: right invertibility of w(r_-n)} that $W(r_{- n})$ is invertible from the right and hence surjective. On the other hand, Lemma \ref{lem: kernel of W(r_-n)} reveals that the kernel of $W(r_{- n})$ is finite-dimensional with
\[
\dim \Ker W(r_{- n}) = n.
\]
Combining all of the above remarks, we arrive at the conclusion that $W(r_{- n})$ is Fredholm on $X(\R_{+})$ with
\[
\Ind W(r_{- n}) = n - 0 = n.
\]
As an immediate consequence of the proved part of the result, we obtain that $W(r_{- n})$ is Fredholm on $X'(\R_{+})$ alongside its Fredholm index formula. Under these conditions, \cite[Chap.~1, Theorem~3.5]{P78} yields that $W(r_{n}) = (W(r_{- n}))^{\ast}$ (see Lemma \ref{lem: adjoint of W(a)}) is also Fredholm on $X''(\R_{+}) = X(\R_{+})$ (see, e.g., \cite[Chap.~1, Theorem~2.7]{BS88}) with 
\[
\Ind W(r_{n}) = - \Ind W(r_{- n}) = - n. 
\]
\end{proof}

%%%%%%%%%%%%%%%%%%%%%%%%%%%%%%%%%%%%%%%%%%%%%%%%%%%%%%%%%%%%%%%%%%%%%%%%%%%%

\subsection{Fundamental homotopy}

The proof of the next result is based on the argument employed in the proof of \cite[Theorem 3.2]{D79}. We provide details for the convenience of the reader.

%%%%%%%%%%%%%%%%%%%%%%%%%%%%%%%%%%%%%%%%%%%%%%%%%%%%%%%%%%%%%%%%%%%%%%%%%%%%

\begin{theorem}\label{thm: fundamental homotopy}
Let $X(\R)$ be a separable rearrangement-invariant Banach function space with nontrivial Boyd indices or a variable Lebesgue space with $p(\cdot)\in \B_{M}(\R)$ and $b\in P\ell \cap C(\dot{\R})$. If the symbol $b$ is elliptic, then $b$ is homotopic to $r_{\wind b}$ in $\M_{X}$.
\end{theorem}

\begin{proof}
Fix $b\in P\ell \cap C(\dot{\R})$ and suppose that $b$ is elliptic. Set $\kappa := \wind b$ and consider the mapping $\mathds{h}(t) := h_{t}$, $t\in [0, 1]$, where each function $h_{t}$ is defined as
\[
h_{t}(\xi) := \big(r_{- \kappa}(\xi) b(\xi)\big)^{t} r_{\kappa}(\xi), 
\qquad 
\xi\in \R.
\]
Note that $\mathds{h}(0) = r_{\kappa}$ and $\mathds{h}(1) = b$. 

\textbf{Claim I:} \textit{$\mathds{h}(t)$ belongs to $BV(\R) \cap C(\dot{\R})$ and is elliptic}. Fix $t\in [0, 1]$. For convenience, denote 
\[
f(\xi) := r_{- \kappa}(\xi) b(\xi), 
\qquad 
g_{t}(\xi) := \xi^{t}. 
\]
Under these conditions, we have $h_{t} = (g_{t} \circ f) r_{\kappa}$. Since $r_{- \kappa}, b\in BV(\R) \cap C(\dot{\R})$, we also have $f\in BV(\R) \cap C(\dot{\R})$. The property $|r_{\pm \kappa}| = 1$ combined with the assumption that $b$ is elliptic implies that $f$ is also elliptic. Taking into account that the functions $f$ and $g_{t} \circ f$ are continuously differentiable up to a finite set of points, it follows from the additivity of the variation and \cite[Chap.~IV, Proposition~1.3]{C78} that
\begin{align*}
V(g_{t} \circ f) = \int_{\R} |(g_{t} \circ f)'(\xi)| \ d\xi &= t \int_{\R} \dfrac{|f'(\xi)|}{|f(\xi)|^{1 - t}} \ d\xi \\
&\le \dfrac{t}{\Big(\displaystyle\inf_{\xi\in \R} |f(\xi)|\Big)^{1 - t}} \int_{\R} |f'(\xi)| \ d\xi = \dfrac{tV(f)}{\Big(\displaystyle\inf_{\xi\in \R} |b(\xi)|\Big)^{1 - t}},
\end{align*}
which in turn reveals that $g_{t} \circ f\in BV(\R)$ and
\begin{equation}\label{eq: estimate on variation}
V((r_{- \kappa}b)^{t}) \le \dfrac{t V(r_{- \kappa}b)}{\Big(\displaystyle\inf_{\xi\in \R} |b(\xi)|\Big)^{1 - t}}.
\end{equation}
Combining this with the continuity of $f$ and $g_{t}$, we get that $g_{t} \circ f\in BV(\R) \cap C(\dot{\R})$ and thus $h_{t}\in BV(\R) \cap C(\dot{\R})$. Similar considerations as the ones made above allow us to conclude that $\mathds{h}(t)$ is elliptic.

\textbf{Claim II:} \textit{$\mathds{h}$ is continuous in $L^{\infty}(\R)$.} Fix $t_{0}\in [0, 1]$. Since the function $r_{- \kappa}b$ is elliptic, for every $t\in [0, 1]$ sufficiently close to $t_{0}$ and for all $\xi\in \R$,
\[
(r_{- \kappa}(\xi) b(\xi))^{|t - t_{0}|} - 1 = \Log(r_{- \kappa}(\xi) b(\xi)) \int_{\gamma([0, 1])} (r_{- \kappa}(\xi) b(\xi))^{z} \ dz,
\]
where $\gamma(x) := x|t - t_{0}|$, $x\in [0, 1]$, and $\Log$ stands for the principal branch of the logarithm. Consequently,
\begin{align*}
|(r_{- \kappa}(\xi) b(\xi))^{|t - t_{0}|} - 1| &\le |\Log(r_{- \kappa}(\xi) b(\xi))| \sup_{z\in \gamma([0, 1])} |(r_{- \kappa}(\xi) b(\xi))^{z}| \, V(\gamma) \\
&= |\ln |b(\xi)| + i\Arg(r_{- \kappa}(\xi) b(\xi))| \sup_{x\in [0, 1]} |b(\xi)|^{x|t - t_{0}|} \, |t - t_{0}| \\
&\le (|\ln |b(\xi)| \, | + |\Arg(r_{- \kappa}(\xi) b(\xi))|) \sup_{x\in [0, 1]} \| b \|_{L^{\infty}(\R)}^{x|t - t_{0}|} \, |t - t_{0}| \\
&\le (|\ln \| b \|_{L^{\infty}(\R)} \, | + \pi) \max\left\{1, \| b \|_{L^{\infty}(\R)}^{|t - t_{0}|}\right\} |t - t_{0}|.
\end{align*}
Therefore,
\begin{equation}\label{eq: limit in Linfty}
\lim_{t \to t_{0}} \| (r_{- \kappa} b)^{|t - t_{0}|} - 1 \|_{L^{\infty}(\R)} = 0.
\end{equation}
On account of this and the estimates
\begin{align*}
\| \mathds{h}(t) - \mathds{h}(t_{0}) \|_{L^{\infty}(\R)} &= \| (r_{- \kappa} b)^{t} r_{\kappa} - (r_{- \kappa} b)^{t_{0}} r_{\kappa} \|_{L^{\infty}(\R)} \\
&\le \| (r_{- \kappa} b)^{\min\{t, t_{0}\}} \|_{L^{\infty}(\R)} \| (r_{- \kappa} b)^{|t - t_{0}|} - 1 \|_{L^{\infty}(\R)} \| r_{\kappa} \|_{L^{\infty}(\R)} \\
&= \| r_{- \kappa} b \|_{L^{\infty}(\R)}^{\min\{t, t_{0}\}} \| (r_{- \kappa} b)^{|t - t_{0}|} - 1 \|_{L^{\infty}(\R)} \\
&= \| b \|_{L^{\infty}(\R)}^{\min\{t, t_{0}\}} \| (r_{- \kappa} b)^{|t - t_{0}|} - 1 \|_{L^{\infty}(\R)},
\end{align*}
a passage to the limit yields that
\[
\lim_{t \to t_{0}} \| \mathds{h}(t) - \mathds{h}(t_{0}) \|_{L^{\infty}(\R)} = 0,
\]
and thus proving that $\mathds{h}$ is continuous in $L^{\infty}(\R)$.

\textbf{Claim III:} \textit{$\mathds{h}$ is continuous in $BV(\R)$.} Fix $t_{0}\in [0, 1]$. Proceeding analogously as in the proof of Claim II, one deduces that
\begin{align*}
\| \mathds{h}(t) - \mathds{h}(t_{0}) \|_{BV(\R)} &= \| (r_{- \kappa} b)^{t} r_{\kappa} - (r_{- \kappa} b)^{t_{0}} r_{\kappa} \|_{BV(\R)} \\
&\le \| (r_{- \kappa} b)^{\min\{t, t_{0}\}} \|_{BV(\R)} \| (r_{- \kappa} b)^{|t - t_{0}|} - 1 \|_{BV(\R)} \| r_{\kappa} \|_{BV(\R)}.
\end{align*}
Employing \eqref{eq: estimate on variation}, we find that the term
\begin{align*}
\| (r_{- \kappa} b)^{\min\{t, t_{0}\}} \|_{BV(\R)} &= \| (r_{- \kappa} b)^{\min\{t, t_{0}\}} \|_{L^{\infty}(\R)} + V((r_{- \kappa} b)^{\min\{t, t_{0}\}}) \\
&\le \| b \|_{L^{\infty}(\R)}^{\min\{t, t_{0}\}} + \dfrac{\min\{t, t_{0}\} V(r_{- \kappa} b)}{\Big(\displaystyle\inf_{\xi\in \R} |b(\xi)|\Big)^{1 - \min\{t, t_{0}\}}} 
\end{align*}
is bounded. On the other hand, \eqref{eq: estimate on variation} also gives that
\[
V((r_{- \kappa} b)^{|t - t_{0}|} - 1) = V((r_{- \kappa} b)^{|t - t_{0}|}) \le \dfrac{|t - t_{0}| V(r_{- \kappa} b)}{\Big(\displaystyle\inf_{\xi\in \R} |b(\xi)|\Big)^{1 - |t - t_{0}|}},
\]
the latter of which when combined with \eqref{eq: limit in Linfty} implies that
\[
\lim_{t \to t_{0}} \| (r_{- \kappa} b)^{|t - t_{0}|} - 1 \|_{BV(\R)} = 0.
\]
Following this sequence of events, the estimate obtained initially leads us to the conclusion that
\[
\lim_{t \to t_{0}} \| \mathds{h}(t) - \mathds{h}(t_{0}) \|_{BV(\R)} = 0,
\]
proving our claim.

Now we are ready to show that $\mathds{h}$ is continuous in $\M_{X}$. Fix $t_{0}\in [0, 1]$. Since by Claim I, we have $\mathds{h}(t)\in BV(\R)$, a direct application of \cite[Theorem~4.3]{K15} gives
\[
\| \mathds{h}(t) - \mathds{h}(t_{0}) \|_{\M_{X}} \le c_{X} \| \mathds{h}(t) - \mathds{h}(t_{0}) \|_{BV(\R)},
\]
which when combined with the result proved in Claim III yields that
\[
\lim_{t \to t_{0}} \| \mathds{h}(t) - \mathds{h}(t_{0}) \|_{\M_{X}} = 0,
\]
as intended.
\end{proof}

%%%%%%%%%%%%%%%%%%%%%%%%%%%%%%%%%%%%%%%%%%%%%%%%%%%%%%%%%%%%%%%%%%%%%%%%%%%%

\subsection{Proof of Theorem \ref{thm: main result}}

\begin{proof}
The proof follows the argument used in \cite[Theorem~3.2]{D79}. Let $a\in C_{X}(\dot{\R})$ and suppose that $a$ is elliptic. Then $a^{- 1}$ exists and, by Theorem \ref{thm: inverse closedness of C_X}, we have $a^{- 1}\in C_{X}(\dot{\R})$. In virtue of this, a direct application of Theorem \ref{thm: compactness of semi-commutators between WH operators with continuous symbols} gives that the operators
\[
I - W(a)W(a^{- 1}), 
\qquad 
I - W(a^{- 1})W(a),
\]
are compact on $X(\R_{+})$ and, consequently, \cite[Chap.~XI, Theorem~5.1]{GGK90} gives that $W(a)$ is a Fredholm operator on $X(\R_{+})$. We now turn our attention to the index formula. Fix $\varepsilon > 0$. On account of Theorems \ref{thm: C_X in terms of piecewise linear - RI BFS version} and \ref{thm: C_p() in terms of piecewise linear}, there is some $b\in P\ell \cap C(\dot{\R})$ such that
\[
\| a - b \|_{\M_{X}} < \varepsilon,
\]
and, consequently,
\[
\| W(a) - W(b) \|_{\B(X(\R_{+}))} < \varepsilon. 
\]
Since $\varepsilon$ is arbitrary and the set of all Fredholm operators is open (see, e.g., \cite[Chap.~XI, Theorem~4.1]{GGK90}), we can choose $\varepsilon$ small enough so that $W(b)$ becomes a Fredholm operator with
\[
\Ind W(a) = \Ind W(b).
\]
On the other hand, it follows from \cite[Theorem~2.3]{KS21} that
\[
\| a - b \|_{L^{\infty}(\R)} < \varepsilon.
\]
In light of this inequality, we can take $\varepsilon$ small enough so that $b$ is elliptic, and $\wind a = \wind b$ (see, e.g., \cite[Lemma~2.7.21]{W12}). On account of the above considerations, we choose $\varepsilon > 0$ sufficiently small so that the following properties hold:
\begin{enumerate}
\item $W(b)$ is a Fredholm operator with $\Ind W(a) = \Ind W(b)$;
\item $b(\xi)\neq 0$ for all $\xi\in \dot{\R}$;
\item $\wind a = \wind b$.
\end{enumerate}
Under these conditions, we are naturally inclined to define
\[
\mathds{H}(t) := W(\mathds{h}(t)), 
\qquad 
t\in [0, 1],
\]
where the symbol $\mathds{h}(t)$ is defined as in the proof of Theorem \ref{thm: fundamental homotopy}. It is clear that $\mathds{H}(0) = W(r_{\wind a})$ and $\mathds{H}(1) = W(b)$.

\textbf{Property I:} $\mathds{H}(t)$ is Fredholm on $X(\R_{+})$. Fix $t\in [0, 1]$. By Theorem \ref{thm: fundamental homotopy}, the symbol $\mathds{h}(t)$ is elliptic, and, consequently, the first part of the proof yields that the Wiener-Hopf operator $W(\mathds{h}(t))$ is Fredholm on $X(\R_{+})$.

\textbf{Property II:} the function $\mathds{H}$ is continuous. For all $t_{0}\in [0, 1]$, we have
\[
\| \mathds{H}(t) - \mathds{H}(t_{0}) \|_{\B(X(\R_{+}))} \le \| \mathds{h}(t) - \mathds{h}(t_{0}) \|_{\M_{X}}.
\]
Faced with this, we are left to pass to the limit as $t\to t_{0}$ and to apply Theorem \ref{thm: fundamental homotopy} to conclude that $\mathds{H}$ is indeed continuous. 

Therefore, $\mathds{H}$ constitutes a homotopy between the Wiener-Hopf operators $W(r_{\wind a})$ and $W(b)$. Finally, since the Fredholm index is homotopically invariant (see, e.g., \cite[Chap.~I, Theorem~3.11]{MP86}), we have
\[
\Ind W(a) = \Ind W(b) = \Ind \mathds{H}(1) = \Ind \mathds{H}(0) = \Ind W(r_{\wind a}) = - \wind a,
\]
where the last equality is justified by Theorem \ref{thm: Fredholmness of W(r_n)}.

Conversely, suppose that $W(a)$ is Fredholm on $X(\R_{+})$ and assume by contradiction that there is some $\xi_{0}\in \dot{\R}$ such that $a(\xi_{0}) = 0$. The remarks made in the first part of the proof allow us to, without loss of generality, consider $a$ to be a continuous piecewise linear function such that
\[
\forall \xi\in \dot{\R}\setminus \{\xi_{0}\} \quad a(\xi) \neq 0.
\]
Under these conditions, we are able to slightly move around the continuous closed curve $a(\dot{\R})$ in such a way that its perturbed versions $a_{\pm \varepsilon}(\dot{\R})$ will either contain the point $0$ in the region enclosed by it or $0$ will lie outside of it. This translates into the existence of some $v\in \C$ (direction) and some small $\varepsilon > 0$ (intensity of movement) such that
\[
a_{\pm \varepsilon}(\xi) := a(\xi) \pm \varepsilon v \neq 0,
\]
for all $\xi\in \dot{\R}$, with
\[
|\wind(a_{- \varepsilon}) - \wind(a_{+ \varepsilon})| = 1.
\]
But then the proved part of the theorem gives us that the operators $W(a_{\pm \varepsilon})$ are Fredholm on $X(\R_{+})$ and satisfy
\[
|\Ind W(a_{+ \varepsilon}) - \Ind W(a_{- \varepsilon})| = 1.
\]
On the other hand,
\[
\| W(a_{+ \varepsilon}) - W(a_{- \varepsilon}) \|_{\B(X(\R_{+}))} = \| W(2\varepsilon v) \|_{\B(X(\R_{+}))} = 2\varepsilon |v|.
\]
Employing the fact that the Fredholm index is stable under small perturbations (see, e.g., \cite[Chap.~XI, Theorem~4.1]{GGK90}), we get
\[
\Ind W(a_{+ \varepsilon}) = \Ind W(a_{- \varepsilon}),
\]
which yields a contradiction. Therefore $a(\xi) \neq 0$ for all $\xi\in \dot{\R}$.
\end{proof}

%%%%%%%%%%%%%%%%%%%%%%%%%%%%%%%%%%%%%%%%%%%%%%%%%%%%%%%%%%%%%%%%%%%%%%%%%%%%
%%%%%%%%%%%%%%%%%%%%%%%%%%%%%%%%%%%%%%%%%%%%%%%%%%%%%%%%%%%%%%%%%%%%%%%%%%%%
%%%%%%%%%%%%%%%%%%%%%%%%%%%%%%%%%%%%%%%%%%%%%%%%%%%%%%%%%%%%%%%%%%%%%%%%%%%%

\section*{Declarations}
\subsection*{Acknowledgements} 

I would like to express my sincere gratitude towards my PhD advisor Oleksiy Karlovych for his attention to the work, constructive criticism and valuable feedback which contributed to enhancing the quality of the paper. I would also like to thank the anonymous referees for their careful reading of the manuscript and their useful remarks.

%%%%%%%%%%%%%%%%%%%%%%%%%%%%%%%%%%%%%%%%%%%%%%%%%%%%%%%%%%%%%%%%%%%%%%%%%%%%

\subsection*{Funding} 

This work is funded by national funds through the FCT – Fundação para a Ciência e a Tecnologia, I.P., under the scope of the projects UIDB/00297/2020 (\url{https://doi.org/10.54499/UIDB/00297/2020}) and UIDP/00297/2020 (\url{https://doi.org/10.54499/UIDP/00297/2020}) (Center for Mathematics and Applications), as well as the PhD scholarship 2022.12247.BD (\url{https://doi.org/10.54499/2022.12247.BD}).

%%%%%%%%%%%%%%%%%%%%%%%%%%%%%%%%%%%%%%%%%%%%%%%%%%%%%%%%%%%%%%%%%%%%%%%%%%%%

\subsection*{Conflict of interest}

The author declares no competing interests.

%%%%%%%%%%%%%%%%%%%%%%%%%%%%%%%%%%%%%%%%%%%%%%%%%%%%%%%%%%%%%%%%%%%%%%%%%%%%

\subsection*{Data availability}

Data sharing not applicable to this article as no datasets were generated or analysed during the current study.

%%%%%%%%%%%%%%%%%%%%%%%%%%%%%%%%%%%%%%%%%%%%%%%%%%%%%%%%%%%%%%%%%%%%%%%%%%%%
%%%%%%%%%%%%%%%%%%%%%%%%%%%%%%%%%%%%%%%%%%%%%%%%%%%%%%%%%%%%%%%%%%%%%%%%%%%%
%%%%%%%%%%%%%%%%%%%%%%%%%%%%%%%%%%%%%%%%%%%%%%%%%%%%%%%%%%%%%%%%%%%%%%%%%%%%

\bibliography{MV1} 

\end{document}